# A Simple Proof for the Four-Color Theorem

Wei-Chang Yeh
Integration & Collaboration Laboratory
Department of Industrial Engineering and Management Engineering
National Tsing Hua University, Hsinchu, Taiwan
yeh@ieee.org

The four-color theorem states that no more than four colors are required to color all nodes in planar graphs such that no two adjacent nodes are of the same color. The theorem was first propounded by Francis Guthrie in 1852. Since then, scholars have either failed to solve this theorem or required computer assistance to prove it. Hence, the goal of this paper is to provide the first correct proof of this 170-year-old mathematical problem composed with the human brain and without computer assistance in only five pages.

**The Triangular Graph**

In graph theory, a planar graph can be drawn on the plane in such a way that no edges cross each other, while a triangular graph is one in which each bounded face is surrounded by three edges. In this paper, the original planar graph must first be transferred into a triangular graph, and each bounded face is called a triangle, including three different nodes and one face in the triangular graphs. Note that the above transformation is not unique. The following lemma explains why each planar graph needs to be transferred into a triangular graph in the paper.

**Lemma 1**   The planar graph is 4-colorable if its related triangular graph is 4-colorable.

**Proof**   Let $G(V, E^*)$ be a planar graph with a node set $V$ and edge set $E^*$ and $G(V, E)$ be the related triangular graph after adding edges into $G(V, E^*)$. If $G(V, E)$ is colored into four different colors, i.e., 4-colorable, such that two adjacent nodes are colored by two different colors, then the removal of any edge from $G(V, E)$ is also 4-colorable.

After removing these edges in $E^*-E$, all nodes colored in $G(V, E)$ will have the same color in $G(V, E^*)$, i.e., all adjacent nodes will be also colored by two different colors in $G(V, E^*)$. Hence, this lemma is true. □

For example, Fig. 1(a) is a planar graph $G(V, E^*)$ and each node is labeled by a number; Fig. 1(b) $G(V, E)$ is a triangular graph formed by adding dashed edges $e_{1,3}$ into $G(V, E^*)$, where $e_{i,j}$ is an undirected arc between the adjacent nodes $i$ and $j$.

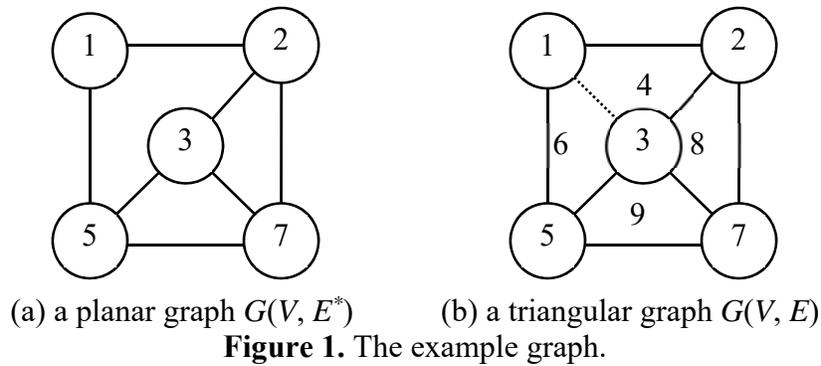

(a) a planar graph $G(V, E^*)$   (b) a triangular graph $G(V, E)$
**Figure 1.** The example graph.

**The Novel Concept of the Complementary Colored Face**

One of the novelties of this paper is to present a brand-new concept called the complementary colored face by letting the color of the face be the color that is not used to color all three nodes in the related triangle. Let $\Delta_{a,b,c,d}$ be a triangle surrounded by nodes $a$, $b$, and $c$ with face $d$ inside, and let $\alpha$, $\beta$, $\chi$, $\delta$ be the four colors. For example, in $\Delta_{1,2,3,4}$ of Fig. 1(b), face (component) 4 must be colored into color $\delta$ if nodes 1, 2, and 3 are colored in $\alpha$, $\beta$, and $\chi$, respectively.

The following lemma discusses an important property of the complementary colored face, which is the key to proving the four-color theorem.

**Lemma 2**   Each color is used once and only once in each triangle if there are only four different colors.

**Proof**   Two adjacent nodes are in different colors, and three different colors are needed to color three different nodes in each triangle. Furthermore, the face must be colored in the color that is different to that of all three nodes. Hence, this lemma is true.   □

**The Novel Four-Component, Four-Color Linear System**

In addition to the novel complementary colored face, another novelty called the four-component, four-color linear system is proposed to transfer the four-color problem into a linear system. This is the reason why the proposed proof can solve the four-color theorem in five pages without the help from computers.

Let zero-one variable $x_{i,j} = 1$ if the $i$th component, i.e., a node or face, is colored in the $j$th color in the triangular graph. Otherwise, let $x_{i,j} = 0$ for $i = 1, 2, \ldots, |V|+|F|-1$ and $j = 1, 2, 3, 4$, where $|V|$ and $|F|$ are the number of nodes and faces, respectively. Hence, we have $4 * (|V|+|F|-1)$ zero-one variables. Because each component can be colored by one and only one color according to the definition of the four-color theorem, and because each color can be used once and only once in each triangle from Lemma 2, we have

$$\sum_{j=1}^{4} x_{i,j} = 1 \text{ for all components } i \text{ in the same triangle} \quad (1)$$

$$\sum_{i} x_{i,j} = 1 \text{ for all components } i \text{ in the same triangle and } j = 1, 2, 3, 4. \quad (2)$$

There are eight linear equations for each triangle. Here, Eq.(1) is called the component equation and Eq.(2) is called the color equation. For example, in Fig. 1(b), all related linear equations are listed in the 2nd column of Table 1, where the first four equations are component equations derived from Eq.(1) and the rest are node equations based on Eq.(2) for each triangle.

For $\sum_{i}\sum_{j=1}^{4} x_{i,j} = \sum_{j=1}^{4}\sum_{i} x_{i,j} = 1$ for all components $i$ in the same triangle and $j = 1, 2, 3, 4$.

One equation (e.g., the last one) is a linear combination of the others for each triangle and can be removed. Afterwards, the last item of all component equations, e.g., $x_{1,4}$, $x_{2,4}$, $x_{3,4}$, $x_{4,4}$, in $\Delta_{1,2,3,4}$ of Table 1, is not in the color equations, and the rest of the equations are linearly independent.

Table 1. The linear equations in Fig. 1.

| Triangle | Component equations(1) and node equations (2) | Last item removed | Duplicates removed |
|---|---|---|---|
| $\Delta_{1,2,3,4}$ | $x_{1,1}+x_{1,2}+x_{1,3}+x_{1,4}=1$<br>$x_{2,1}+x_{2,2}+x_{2,3}+x_{2,4}=1$<br>$x_{3,1}+x_{3,2}+x_{3,3}+x_{3,4}=1$<br>$x_{4,1}+x_{4,2}+x_{4,3}+x_{4,4}=1$<br>$x_{1,1}+x_{2,1}+x_{3,1}+x_{4,1}=1$<br>$x_{1,2}+x_{2,2}+x_{3,2}+x_{4,2}=1$<br>$x_{1,3}+x_{2,3}+x_{3,3}+x_{4,3}=1$<br>$x_{1,4}+x_{2,4}+x_{3,4}+x_{4,4}=1$ | $x_{1,1}+x_{1,2}+x_{1,3}+x_{1,4}=1$<br>$x_{2,1}+x_{2,2}+x_{2,3}+x_{2,4}=1$<br>$x_{3,1}+x_{3,2}+x_{3,3}+x_{3,4}=1$<br>$x_{4,1}+x_{4,2}+x_{4,3}+x_{4,4}=1$<br>$x_{1,1}+x_{2,1}+x_{3,1}+x_{4,1}=1$<br>$x_{1,2}+x_{2,2}+x_{3,2}+x_{4,2}=1$<br>$x_{1,3}+x_{2,3}+x_{3,3}+x_{4,3}=1$ | $x_{1,1}+x_{1,2}+x_{1,3}+x_{1,4}=1$<br>$x_{2,1}+x_{2,2}+x_{2,3}+x_{2,4}=1$<br>$x_{3,1}+x_{3,2}+x_{3,3}+x_{3,4}=1$<br>$x_{4,1}+x_{4,2}+x_{4,3}+x_{4,4}=1$<br>$x_{1,1}+x_{2,1}+x_{3,1}+x_{4,1}=1$<br>$x_{1,2}+x_{2,2}+x_{3,2}+x_{4,2}=1$<br>$x_{1,3}+x_{2,3}+x_{3,3}+x_{4,3}=1$ |
| $\Delta_{1,3,5,6}$ | $x_{1,1}+x_{1,2}+x_{1,3}+x_{1,4}=1$<br>$x_{3,1}+x_{3,2}+x_{3,3}+x_{3,4}=1$<br>$x_{5,1}+x_{5,2}+x_{5,3}+x_{5,4}=1$<br>$x_{6,1}+x_{6,2}+x_{6,3}+x_{6,4}=1$<br>$x_{1,1}+x_{3,1}+x_{5,1}+x_{6,1}=1$<br>$x_{1,2}+x_{3,2}+x_{5,2}+x_{6,2}=1$<br>$x_{1,3}+x_{3,3}+x_{5,3}+x_{6,3}=1$<br>$x_{1,4}+x_{3,4}+x_{5,4}+x_{6,4}=1$ | $x_{1,1}+x_{1,2}+x_{1,3}+x_{1,4}=1$<br>$x_{3,1}+x_{3,2}+x_{3,3}+x_{3,4}=1$<br>$x_{5,1}+x_{5,2}+x_{5,3}+x_{5,4}=1$<br>$x_{6,1}+x_{6,2}+x_{6,3}+x_{6,4}=1$<br>$x_{1,1}+x_{3,1}+x_{5,1}+x_{6,1}=1$<br>$x_{1,2}+x_{3,2}+x_{5,2}+x_{6,2}=1$<br>$x_{1,3}+x_{3,3}+x_{5,3}+x_{6,3}=1$ | $x_{5,1}+x_{5,2}+x_{5,3}+x_{5,4}=1$<br>$x_{6,1}+x_{6,2}+x_{6,3}+x_{6,4}=1$<br>$x_{1,1}+x_{3,1}+x_{5,1}+x_{6,1}=1$<br>$x_{1,2}+x_{3,2}+x_{5,2}+x_{6,2}=1$<br>$x_{1,3}+x_{3,3}+x_{5,3}+x_{6,3}=1$ |
| $\Delta_{2,3,7,8}$ | $x_{2,1}+x_{2,2}+x_{2,3}+x_{2,4}=1$<br>$x_{3,1}+x_{3,2}+x_{3,3}+x_{3,4}=1$<br>$x_{7,1}+x_{7,2}+x_{7,3}+x_{7,4}=1$<br>$x_{8,1}+x_{8,2}+x_{8,3}+x_{8,4}=1$<br>$x_{2,1}+x_{3,1}+x_{7,1}+x_{8,1}=1$<br>$x_{2,2}+x_{3,2}+x_{7,2}+x_{8,2}=1$<br>$x_{2,3}+x_{3,3}+x_{7,3}+x_{8,3}=1$<br>$x_{2,4}+x_{3,4}+x_{7,4}+x_{8,4}=1$ | $x_{2,1}+x_{2,2}+x_{2,3}+x_{2,4}=1$<br>$x_{3,1}+x_{3,2}+x_{3,3}+x_{3,4}=1$<br>$x_{7,1}+x_{7,2}+x_{7,3}+x_{7,4}=1$<br>$x_{8,1}+x_{8,2}+x_{8,3}+x_{8,4}=1$<br>$x_{2,1}+x_{3,1}+x_{7,1}+x_{8,1}=1$<br>$x_{2,2}+x_{3,2}+x_{7,2}+x_{8,2}=1$<br>$x_{2,3}+x_{3,3}+x_{7,3}+x_{8,3}=1$ | $x_{7,1}+x_{7,2}+x_{7,3}+x_{7,4}=1$<br>$x_{8,1}+x_{8,2}+x_{8,3}+x_{8,4}=1$<br>$x_{2,1}+x_{3,1}+x_{7,1}+x_{8,1}=1$<br>$x_{2,2}+x_{3,2}+x_{7,2}+x_{8,2}=1$<br>$x_{2,3}+x_{3,3}+x_{7,3}+x_{8,3}=1$ |
| $\Delta_{3,5,7,9}$ | $x_{3,1}+x_{3,2}+x_{3,3}+x_{3,4}=1$<br>$x_{5,1}+x_{5,2}+x_{5,3}+x_{5,4}=1$<br>$x_{7,1}+x_{7,2}+x_{7,3}+x_{7,4}=1$<br>$x_{9,1}+x_{9,2}+x_{9,3}+x_{9,4}=1$<br>$x_{3,1}+x_{5,1}+x_{7,1}+x_{9,1}=1$<br>$x_{3,2}+x_{5,2}+x_{7,2}+x_{9,2}=1$<br>$x_{3,3}+x_{5,3}+x_{7,3}+x_{9,3}=1$<br>$x_{3,4}+x_{5,4}+x_{7,4}+x_{9,4}=1$ | $x_{3,1}+x_{3,2}+x_{3,3}+x_{3,4}=1$<br>$x_{5,1}+x_{5,2}+x_{5,3}+x_{5,4}=1$<br>$x_{7,1}+x_{7,2}+x_{7,3}+x_{7,4}=1$<br>$x_{9,1}+x_{9,2}+x_{9,3}+x_{9,4}=1$<br>$x_{3,1}+x_{5,1}+x_{7,1}+x_{9,1}=1$<br>$x_{3,2}+x_{5,2}+x_{7,2}+x_{9,2}=1$<br>$x_{3,3}+x_{5,3}+x_{7,3}+x_{9,3}=1$ | $x_{9,1}+x_{9,2}+x_{9,3}+x_{9,4}=1$<br>$x_{3,1}+x_{5,1}+x_{7,1}+x_{9,1}=1$<br>$x_{3,2}+x_{5,2}+x_{7,2}+x_{9,2}=1$<br>$x_{3,3}+x_{5,3}+x_{7,3}+x_{9,3}=1$ |

Additionally, duplicate component equations can be removed. E.g., $x_{3,1}+x_{3,2}+x_{3,3}+x_{3,4}=1$ appears four times: $\Delta_{1,2,3,4}$, $\Delta_{1,3,5,6}$, $\Delta_{2,3,7,8}$, and $\Delta_{3,5,7,9}$ in the 3rd column of Table 1. After that, the last item of each linear equation never appears in any equation for other triangles. In other words, all linear equations are linearly independent in the linear system, e.g., $x_{1,4}$, $x_{2,4}$, $x_{3,4}$, $x_{4,4}$, $x_{4,1}$, $x_{4,2}$ and $x_{4,3}$ only appear in $\Delta_{1,2,3,4}$. The above discussions prove the following lemma.

**Lemma 3**   All row vectors are linear independence of both the coefficient matrix and augmented matrix in the four-component, four-color linear system constructed by the Eqs.(1) and (2) after removing one color equation from each triangle and all duplicates.

**The Proof**

The proof of the four-color theorem is provided below.

**Theorem 1**   Each planar graph is 4-colorable.

**Proof**   Each planar graph can be transferred into a triangular graph by adding required edges, and said triangular graph can be further transferred into a four-component, four-color linear system after introducing the zero-one variable $x_{i,j}$ for $i$ = 1, 2, …, $|V|+|F|-1$ and $j$ = 1, 2, 3, 4. All row vectors are linearly independent in both the coefficient matrix and augmented matrix in the four-component, four-color linear system. Furthermore, all constants are one in all equations in the linear system. Hence, the linear system has integer-feasible solutions, i.e., the triangular graph is four-colorable. Because the removal of edges that are not in the original planar graph from the triangular graph has no effect on the colors in nodes, the four-color theorem is correct.   □

**Discussion**

The four-color theorem is the most actively worked-on problem in graph theory. The first correct, simple, and unassisted proof is provided for the 170-year-old four-color theorem in this paper based on two innovations: the complementary colored face and the four-component, four-color linear system. These novel concepts build a bridge between graph theory and linear algebra that can be harnessed by scholars to assist them in solving other problems.

## References

[1] https://en.wikipedia.org/wiki/Four_color_theorem

[2] Trudeau, Richard J. (1993). Introduction to Graph Theory (Corrected, enlarged republication. ed.). New York: Dover Pub. p. 64. ISBN 978-0-486-67870-2.

[3] L. Kronecker, "Vorlesungen über die Theorie der Determinanten", Leipzig (1903)

[4] A. Capelli, "Sopra la compatibilitá o incompatibilitá di più equazioni di primo grado fra picì incognite" Revista di Matematica, 2 (1892) pp. 54–58